\documentclass{amsart}

\usepackage{amssymb}

\newtheorem{thm}{Theorem}%[section]

\newtheorem{lem}[thm]{Lemma}
\newtheorem{cor}[thm]{Corollary}

\newtheorem{prop}[thm]{Proposition}

\theoremstyle{definition}

\newtheorem{say}[thm]{}
\newtheorem{exmp}[thm]{Example}

\newtheorem{rem}[thm]{Remark}          
\newtheorem{ack}{Acknowledgments}

\newtheorem{defn-thm}[thm]{Definition--Theorem}  %!!!!!!!!!!!!!!!!!!!!!!!!
\theoremstyle{remark}

\renewcommand{\c}[0]{{\mathbb C}}  

\renewcommand{\o}[0]{{\mathcal O}} 
\newcommand{\z}[0]{{\mathbb Z}}

\newcommand{\p}[0]{{\mathbb P}}
\newcommand{\f}[0]{{\mathbb F}}
\newcommand{\q}[0]{{\mathbb Q}}

\newcommand{\qtq}[1]{\quad\mbox{#1}\quad}

\newcommand{\pic}[0]{\operatorname{Pic}}

\newcommand{\rank}[0]{\operatorname{rank}}
\newcommand{\mult}[0]{\operatorname{mult}}

\newcommand{\supp}[0]{\operatorname{Supp}}

\newcommand{\im}[0]{\operatorname{im}}

\newcommand{\aut}[0]{\operatorname{Aut}}

\newcommand{\sing}[0]{\operatorname{Sing}}

\newcommand{\tors}[0]{\operatorname{tors}}

\newcommand{\weil}[0]{\operatorname{Weil}}

\begin{document}
\bibliographystyle{amsalpha}

\title{Positive Sasakian structures on  5-manifolds}
\author{J\'anos Koll\'ar}
%\today
\maketitle

A quasi--regular Sasakian structure on a manifold $L$
is equivalent to writing $L$ as the unit circle subbundle
of a holomorphic Seifert $\c^*$-bundle over a complex algebraic
orbifold $(X,\Delta=\sum (1-\frac1{m_i})D_i)$. The Sasakian structure is called
positive if the orbifold first Chern class 
$c_1(X)-\Delta=-(K_X+\Delta)$ is positive. We are especially interested
in the case when the Riemannian metric part of the 
Sasakian structure is Einstein. By a result of
\cite{kob, bg00}, this happens iff
\begin{enumerate}
\item $-(K_X+\Delta)$ is positive, 
\item the first Chern class of the Seifert bundle
$c_1(L/X)$ is a rational multiple of $-(K_X+\Delta)$, and
\item there is an orbifold K\"ahler--Einstein metric on $(X,\Delta)$.
\end{enumerate}
If the first two conditions  hold, we say that
$f:L\to (X,\Delta)$ is  a {\it pre-SE} 
(or pre--Sasakian--Einstein) Seifert bundle.

Note that if $H_2(L,\q)=0$ then $H_2(S,\q)\cong \q$
and so  the second condition is automatic.

While it is not true that all pre-SE Seifert bundles
 carry a 
Sasakian--Einstein  structure, all known topological obstructions
to the existence of a Sasakian--Einstein  structure
in dimension 5 are  consequences of the pre-SE condition.

A 2--dimensional orbifold $(S,\Delta)$
such that $-(K_S+\Delta)$ is positive is also called
a {\it log Del Pezzo surface}. For any such, $S$ is a rational surface
with quotient singularities. 
By \cite[2.4]{em5}, Seifert $\c^*$-bundles
over $(S, \sum (1-\frac1{m_i})D_i)$ are
uniquely classified by a homology class
$B\in H_2(S,\z)$ and integers
$0< b_i<m_i$ with $(b_i,m_i)=1$.
We denote the corresponding Seifert  $\c^*$-bundle 
(resp.\ $S^1$-bundle) by
$$
Y\bigl(S,B, \textstyle{\sum} \tfrac{b_i}{m_i}D_i\bigr)
\qtq{resp.} L\bigl(S, B, \textstyle{\sum} \tfrac{b_i}{m_i}D_i\bigr).
$$

\begin{say}[Classification problems] Following
\cite{em5}, we want to address 3 problems.
\begin{enumerate}
\item Given $(S,\Delta)$ and $B$, determine the 5--manifold
 $L(S,B,\sum \tfrac{b_i}{m_i}D_i)$.
\item Given a 5--manifold $L$, decide if it can be
written as $L(S,B,\sum \tfrac{b_i}{m_i}D_i)$ for some  $(S,\Delta)$ and $B$.
\item Given a 5--manifold $L$, describe all of its representations
 as $L(S,B,\sum \tfrac{b_i}{m_i}D_i)$ for some  $(S,\Delta)$ and $B$.
\end{enumerate}

The case when $L$ is simply connected was studied in
\cite{em5}. 
Problem (1) was solved in \cite[5.7]{em5}.
As to Problem (2), note that
by a result of \cite{smale}, such  manifolds are
determined by the second homology group $H_2(L,\z)$
if $w_2=0$. (See \cite{barden} for the $w_2\neq 0$ case.)
The torsion part of $H_2(L,\z)$ can be completely described 
as follows.
(There are partial results about the rank of the free part.)

\begin{thm}\label{h2.tors.thm} \cite[1.4]{em5}
Let $L$ be a compact 5--manifold  
such that $H_1(L,\z)=0$.
Assume that $L$ has a positive Sasakian structure
$f:L\to (S,\Delta)$. 
Then  the torsion subgroup of the second homology, $\tors H_2(L,\z)$,
is determined by $(S,\Delta)$ and it
is one of the following.
\begin{enumerate}
\item $(\z/m)^2$ for some $m\geq 1$,
\item $(\z/5)^4$ or $(\z/4)^4$,
\item $(\z/3)^4,(\z/3)^6$ or $(\z/3)^8$,
\item $(\z/2)^{2n}$ for some $n\geq 1$.
\end{enumerate}
\end{thm}

Note that the theorem gives restrictions on $L$ but does not say anything about
$(S,\Delta)$. 
While it  seems hopelessly complicated to describe all
Seifert bundle  structures on $S^5$, there  are 
few positive Sasakian structures on  complicated
5--manifolds. In particular, Problem (3) was
settled for 
 $\tors H_2(L,\z)\cong (\z/5)^4$  \cite[1.8.2]{em5}.
\end{say}

The aim of this note is fourfold. First, 
 we extend several of the results to
 5--manifolds which are not simply connected. 
Second, I would like to correct
the mistake in the classification of the case 
(\ref{h2.tors.thm}.1) which was brought to my attention by
K.\ Galicki. Third, Problem (3) is solved 
in  the  $(\z/4)^4$ and    $(\z/3)^8$  cases.
Finally, we show that in in many cases, all positive Sasakian
structures are given by hypersurfaces in $\c^4$
and  also give explicit equations
for these. 

Toward the first goal, we show that the groups 
listed in (\ref{h2.tors.thm}.1--4) give large subgroups  of all possible
$\tors H_2(L,\z)$.

\begin{thm} \label{h2.tors.new.thm}
Let $L$ be a compact 5--manifold with a positive
 Sasakian structure. Then $\tors H_2(L,\z)$ has a subgroup
$G$ as in (\ref{h2.tors.thm}.1--4) such that the quotient group
 $\tors H_2(L,\z)/G$ is
  $\z/r$, $\z/r+\z/2$ or $\z/r+\z/3$  for some $r$.
\end{thm}

Our next aim is to prove that in most cases
the existence of the subgroup $G$ controls 
$L$ and the Sasakian structure very tightly.

\begin{thm}\label{h2.tors.gen.cases} 
Let $f:L\to (S,\Delta)$ be a compact 5--manifold with a positive
 Sasakian structure.
 Assume that  $H_2(L,\z)\supset (\z/m)^2$ 
where  $m\geq  28$ or $m\geq 12$ and $(m,6)=1$.
Then 
\begin{enumerate}
\item $S$ is a Del Pezzo surface with cyclic Du Val singularities only
and $\Delta=(1-\frac1{m'})C$ where $C\in |-K_S|$ is a smooth elliptic curve
and $m'$ is a multiple of $m$.
\item there are 132 families of such $(S,C)$ up to deformations that
preserve the singularities,
\item for every such $S$,
$\pi_1(S^0)$ is abelian of order $\leq 9$,
where $S^0:=S\setminus \sing S$ denotes the set of smooth points.
$\pi_1^{orb}(S,(1-\frac1{m'})C)$ is also abelian of order $\leq 9$,
save 3 cases which are nonabelian of order $2^3, 2^4$ or $3^3$.
The precise values are in Table 2.
\item The fundamental group of $L$ is given by  a central extension
$$
0\to \z/d\to \pi_1(L)\to \pi_1^{orb}(S,(1-\tfrac1{m'})C) \to 1,
$$
for some $d$.
\item The torsion of the second homology group of $L$  
sits in an exact sequence
$$
0\to (\z/m')^2\to \tors H_2(L,\z)\to  \pi_1(S^0).
$$
\end{enumerate}
\end{thm}

It is quite likely that the above theorem holds for all
$m\geq 12$. However, for $m<12$, there are many counter examples,
see \cite[6.6--7]{em5}.

\begin{say}[Sasakian--Einstein structures]  In each of the 132 cases the
existence of positive Sasakian and of pre--Sasakian--Einstein structures
over $(S,\Delta)$ is  effectively decidable, but some of the
computations may be lengthy.

There are 93 cases of $(S,C)$ when $\pi_1(S^0)=1$. 
In \cite[1.8.1]{em5} I claimed, incorrectly, that they all
admit a Seifert bundle with Sasakian--Einstein structure.
We see in (\ref{corr181.say})
that these all admit Seifert bundles
with a positive Sasakian  structure,  but
only 19 of them admit a Seifert bundle with pre--SE structure.

For the remaining 39 cases 
with $\pi_1(S^0)\neq 1$
my computations are not yet complete.
There are some cases that do not have any smooth Seifert bundles
over them (\ref{sm.sb.ex.say}), and in many other cases
Sasakian--Einstein structures exist.
\end{say}

\begin{say}[Hypersurface examples]
Concrete examples of Sasakian structures
can be obtained using $\c^*$-actions
on algebraic hypersurfaces.
A linear $\c^*$-action on $\c^n$ can be diagonalized
and so given by weights $(w_1,\dots,w_n)$.
We consider only actions which are effective  (thus the weights
are relatively prime)
and  the origin is an attracting fixed point
(thus $w_i> 0$ for every $i$).
Let $Y\subset \c^n$ be an algebraic variety, smooth away from the origin,
which is invariant under the $\c^*$-action.
Then its link
$L:=Y\cap S^{2n-1}(1)$ inherits a natural quasi--regular Sasakian
structure and every quasi--regular Sasakian
structure arises this way.

This description is especially simple when $Y$ is a hypersurface.
In this case everything is described
by the weights $(w_1,\dots,w_n)$ and a weighted homogeneous
polynomial $f(x_1,\dots,x_n)$ which defines $Y$.

Thus, up to deformations, we need to specify
the weights $(w_1,\dots,w_n)$ and the weighted degree
$d$ of $f$. 
I write
$L^*(w_1,\dots,w_n;d)$ for any such quasi--regular Sasakian
manifold. (The $^*$ is to remind one that using weights
is in some sense dual to the notation in \cite{bgk}.)

The simplest examples where $\pi_1(L)\neq 1$
are created by taking a quotient by a subgroup
of  $\c^*$. 
These are all cyclic. Thus the symbol
$$
L^*(w_1,\dots,w_n;d)/(\z/m)\qtq{or}
L^*(w_1,\dots,w_n;d)/\tfrac1{m}(w_1,\dots,w_n)
$$
stands for any Sasakian manifold obtained
as $L/(\z/m)\subset Y/(\z/m)$ where
$Y:=(f=0)$ is the zero set of a degree $d$ weighted homogeneous
polynomial such that $Y$ is smooth outside the origin
and we take the quotient by the $\z/m$-action generated by
$$
(x_1,\dots,x_n)\mapsto (\epsilon^{w_1}x_1,\dots,\epsilon^{w_n}x_n)
\qtq{where} \epsilon=e^{2\pi i/m}.
$$
We call these the {\it obvious quotients}.

One can also take quotients by groups that are
not contained in $\c^*$. For instance, one can take
quotients
$$
L^*(w_1,\dots,w_n;d)/\bigl(\tfrac1{r}(a_1,\dots,a_n)\times
\tfrac1{m}(w_1,\dots,w_n)\bigr)
$$
where the first factor of $\z/r\times \z/m$
acts via
$$
(x_1,\dots,x_n)\mapsto (\eta^{a_1}x_1,\dots,\eta^{a_n}x_n)
\qtq{where} \eta=e^{2\pi i/r},
$$
and the second factor acts as above.
One should keep in mind that not all hypersurfaces
$L^*(w_1,\dots,w_n;d)$ admit such a $\z/r$-action.
In using this notation, we assume that we consider a case
when the action exists.
\end{say}

We are now ready to state that classification theorem
for those cases when $H_2(L,\z)$ contains a large torsion subgroup.

\begin{thm} \label{main.thm} 
Let $L$ be a compact 5--manifold with a positive
 Sasakian
structure.  
\begin{enumerate}
\item  If $H_2(L,\z)\supset (\z/5)^4$ then 
$H_2(L,\z)= (\z/5)^4$ and 
$L=L^*(5,5,15,6;30)$ or an obvious quotient by $\z/m$ for  $(m,15)=1$.
 The simplest equation is
$$
(x^{6}+y^{6}+z^2+t^{5}=0)/\tfrac1{m}(5,5,15,6).
$$
The moduli space of Sasakian--Einstein structures is naturally
parametrized by the moduli space of  genus 2 curves.

\item  If $H_2(L,\z)\supset (\z/3)^8$ then 
 $H_2(L,\z)=(\z/3)^8$ and 
$L=L^*(3,3,15,10;30)$ or an obvious quotient by $\z/m$ for  $(m,15)=1$.
 The simplest equation is
$$
(x^{10}+y^{10}+z^2+t^{3}=0)/\tfrac1{m}(3,3,15,10).
$$
 The moduli space  of Sasakian--Einstein structures  is naturally
parametrized by the moduli space of  hyperelliptic curves of  genus 4.

\item If $H_2(L,\z)\supset (\z/4)^4$ then 
either $H_2(L,\z)= (\z/4)^4$ or  $H_2(L,\z)= \z+(\z/4)^4$.
In the first case, there are 5  families: 
\begin{enumerate}
\item $L=L^*(4,4,8,5;20)$ or an obvious quotient
by $\z/m$ for $(m,2)=1$,  
 with sample equation 
$$
(x^{5}+y^{5}+yz^2+t^{4}=0)/\tfrac1{m}(4,4,8,5),
$$
\item $Y=L^*(2,2,6,3;12)/\frac12(1,1,1,e)$
where $e\in \{0,1\}$ or an obvious quotient by $\z/m$ for  $(m,6)=1$.
The simplest equation is
$$
(x^6+y^6+z^2+t^4=0)/\bigl(\tfrac12(1,1,1,e)\times \tfrac1{m}(2,2,6,3)\bigr).
$$ 

\item 
 $L=L^*(4,8,20,9;40)$ or an obvious quotient by $\z/m$ for  $(m,6)=1$.
The simplest equation is
$$
(x^{10}+y^{5}+z^2+xt^{4}=0)/\tfrac1{m}(4,8,20,9).
$$
\item $L=L^*(2,4,10,5;20)/\frac12(1,0,1,1)$
or an obvious quotient by $\z/m$ for  $(m,10)=1$.
The simplest equation is
$$
(x^{10}+y^{5}+z^2+t^4=0)/
\bigl(\tfrac12(1,0,1,1)\times \tfrac1{m}(2,4,10,5)\bigr).
$$
\end{enumerate}

In the second case,
there are infinitely many
positive Sasakian structures  with the same $\pi_1$ but only one
which is pre--SE. These are:
\begin{enumerate}\setcounter{enumii}{4}
\item  
$L=L^*(4,4,12,5;24)$ 
or an obvious quotient by $\z/m$ for  $(m,10)=1$.
The simplest equation is
$$
(x^{6}+y^{6}+z^2+xt^4=0)/\tfrac1{m}(4,4,12,5).
$$
\end{enumerate}
The moduli space  of Sasakian--Einstein structures  is naturally
parametrized by the moduli space of  pairs
$(C,p)$ where $\pi:C\to \p^1$ is a genus 2 curve and $p\in C$ is not a
branch point in cases (a),(b),(e) and  $p\in C$ is  a
branch point in cases (c),(d).
\end{enumerate}
\end{thm}

The above parametrizations of the moduli spaces of 
Sasakian--Einstein structures are one--to--one.
However, if we  consider the moduli spaces of
the resulting Einstein metrics, we get a two--to--one
parametrization since a curve and its complex conjugate
give the same Einstein metric, see \cite[18--21]{bgk}.

I see no a priori reason why all these cases should be
realizable as hypersurface quotients. K.\ Galicki  told
 me that he found the
example $(x^{6}+y^{6}+z^2+t^5=0)$
 with  $H_2(L,\z)=(\z/5)^4$.
This lead me to realize that in fact all  examples with 
$H_2(L,\z)=(\z/5)^4$
are hypersurfaces. The equations are also easy to see
in the $(\z/3)^8$ case, but are rather mysterious
from the point of view of my proof for several
of the $(\z/4)^4$ cases.

\section{Reduction to algebraic geometry}

\begin{say}\label{algred.say} As in  \cite{bgk, em5},
to $L=Y\cap S^{2n-1}(1)$ we associate a projective
algebraic variety $X=Y\setminus\{0\}/\c^*$
with cyclic quotient singularities. There is also
a natural $\q$-divisor $\Delta=\sum (1-\tfrac1{m_i})D_i$
on $X$ where $D_i\subset X$ is a divisor such that
the stabilizer of the  $\c^*$-action has order
$m_i$ over the points of  $D_i$.
The positivity of the Sasakian structure is
equivalent to the log Fano condition:
$$
-\bigl(K_X+\sum (1-\tfrac1{m_i})D_i\bigr)\qtq{is ample.}
$$
Thus if $\dim L=5$, we are looking for
pairs $(S,\Delta)$
where 
\begin{enumerate}
\item $S$ is a projective surface with cyclic quotient singularities,
\item $\Delta=\sum (1-\tfrac1{m_i})D_i$ is a $\q$-divisor, and
\item $-(K_S+\Delta)$ is ample.
\end{enumerate}
There are a few more conditions which we do not need for now.

A rather easy result \cite[6.3]{em5} on log Del Pezzo surfaces 
gives the following restriction on the $D_i$:
\begin{enumerate}\setcounter{enumi}{3}
\item there is at most one index $i$ such that $g(D_i)>0$.
\end{enumerate}
Because of the special role of $D_i$, we 
frequently set $D:=D_i$, $r:=m_i$ and 
use the notation $(S,(1-\tfrac1{r})D+\Delta)$.

It turns out that the ampleness of $-(K_S+(1-\tfrac1{r})D+\Delta)$ forces
$g(D)$ or $r$ to be quite small and this accounts for
 the restrictions in Theorem \ref{h2.tors.thm}.

This is relevant to our purposes since
\cite[5.7]{em5} 
 relates the torsion of $H_2(L,\z)$
to the divisor $\sum (1-\tfrac1{m_i})D_i$:

\begin{enumerate}\setcounter{enumi}{4}
\item If $H_1(L,\z)=0$ then $\tors H_2(L,\z)=\sum (\z/m_i)^{2g(D_i)}$.
Thus, in the log Del Pezzo case,  $\tors H_2(L,\z)\cong (\z/r)^{2g(D)}$.
\end{enumerate}

If $H_1(L,\z)\neq 0$, the spectral sequence  in \cite[5.10]{em5}
becomes quite messy, but in the  log Del Pezzo case only
one nonzero map involves $H^3(L,\z)$.
The same proof gives the following weaker result
when the  first ordinary homology of the base is trivial.
\end{say}

\begin{prop} \label{thm2.pf} 
Let $f:L^5\to (S,\Delta=\sum (1-\frac1{m_i})D_i)$
be a smooth Seifert bundle over  a projective surface with
quotient singularities. Assume that $H_1(S,\z)=0$ and
let $s=\rank H^2(S,\q)$. Then
there is an exact sequence
$$
H_1(S^0,\z)\cong H^3(S,\z)\to H^3(L,\z)\to \z^{s-1}+
\sum (\z/m_i)^{2g(D_i)}\to 0. 
\eqno{(\ref{thm2.pf}.1)}
$$
\end{prop}

Proof.
 The argument very closely follows
\cite[5.9--10]{em5}.  
As there, we have exact sequences
$$
0\to R^1f_*\z_L\stackrel{\tau}{\to} \z_S \to Q\to 0
\eqno{(\ref{thm2.pf}.2)}
$$
and 
$$
0\to \sum_i \z_{P_j}/n_j \to Q \to \sum_i \z_{D_i}/m_i\to 0,
\eqno{(\ref{thm2.pf}.3)}
$$
where $P_j\in S$ are the singular points.
This implies that
$H^i(S,Q)=\sum_i H^i(D_i,\z/m_i)$ for $i\geq 1$.
The key piece of the long cohomology sequence of (\ref{thm2.pf}.2)
 is
$$
H^1(S,\z)\to \sum H^1(D_i,\z/m_i)\to H^2(S,R^1f_*\z_L)
\to H^2(S,\z) \to \sum H^2(D_i,\z/m_i)
$$
Here $H^1(S,\z)=0$ and $H^2(S,\z)\cong \z^s$
by assumption.
The right hand group is torsion, hence
there is a noncanonical isomorphism
$$
H^2(S,R^1f_*\z_L)\cong \z^s+\sum_i H^1(D_i,\z/m_i).
$$
Therefore, in
  the Leray spectral sequence
$H^i(S,R^jf_*\z_L)\Rightarrow H^{i+j}(L,\z)$ the $E_2$ term is
$$
\begin{array}{lccll}
\z \quad & (\mbox{torsion}) & \quad\z^{s}+\sum_i (\z/m_i)^{2g(D_i)}\quad & 
(\mbox{torsion})\quad & \z\\
\z & 0 & \z^{s} & H^3(S,\z) & \z.
\end{array}
\eqno{(\ref{thm2.pf}.4)}
$$
Since $H^3(S,\z)\cong H_1(S^0,\z)$ by \cite[4.2]{em5},
we get the exact sequence
$$
H_1(S^0,\z)\to \tors H^3(L,\z)\to \sum_i H^1(D_i,\z/m_i)\to 0. \qed
$$

\begin{cor} Notation and assumptions as in (\ref{thm2.pf}).
Let $p_j\in S$ be the singular points and $n_j$ the order of the local
fundamental group.
\begin{enumerate}
\item If any two of the numbers $m_i, n_j$ are relatively prime
then the sequence (\ref{thm2.pf}.1) is left exact.
\item If, in addition  $|H_1(S^0,\z)|$ is relatively prime to $\prod m_i$
then the sequence (\ref{thm2.pf}.1)  splits.
\end{enumerate}
\end{cor}

Proof. If any two of the numbers $m_i, n_j$ are relatively prime,
 the sequence  in (\ref{thm2.pf}.2)
splits
and $H^0(S,Q)\cong \z/(\prod n_j\cdot \prod m_i)$.
This implies that $H^1(S,R^1f_*\z_L)=0$.
The $E_2$ term of the Leray spectral sequence
(\ref{thm2.pf}.4) is 
$$
\begin{array}{lccll}
\z \quad & 0 & \quad\z^{s}+\sum_i (\z/m_i)^{2g(D_i)}\quad & 
(\mbox{torsion})\quad & \z\\
\z & 0 & \z^{s} & H^3(S,\z) & \z.
\end{array}
$$
Thus $H^3(L,\z)$ sits in an exact sequence
$$
0\to H^3(S,\z) \to H^3(L,\z)\to \z^{s}+\sum_i (\z/m_i)^{2g(D_i)}
\to \z
$$
As before,  $H^3(S,\z)\cong H_1(S^0,\z)$
and the extension of the torsion part splits
if $|H_1(S^0,\z)|$ is relatively prime to $\prod m_i$.\qed
\medskip

This result turns out to be quite useful, since
the groups $H_1(S^0,\z)$ are rather special 
 when
 $(S,\Delta)$ is a log Del Pezzo surface. Thus $S$ is
a rational surface and so $H_1(S,\z)=0$.
 A project of C.\ Xu aims  to give a complete determination of
the possible fundamental groups $\pi_1(S^0)$
where $S$ is a rational surface with quotient singularities.
While this is not easy, the first homology
  $H_1(S^0,\z)$ is not hard to control.

The proof of (\ref{h2.tors.new.thm}) is obtained by combining
(\ref{thm2.pf}) and (\ref{fundgr.lem}).

\begin{lem} \label{fundgr.lem}
Let  $S$ be a log Del Pezzo surface with quotient singularities.
Then  $H_1(S^0,\z)$
is either $\z/m$,  $\z/m+\z/2$ or $\z/m+\z/3$ for some $m\geq 1$.
\end{lem}

Proof. Let $S'\to S$ be the corresponding Galois cover
with  Galois group $G$. The stabilizers of points
are subgroups of the local fundamental groups at the
singularities, and hence cyclic.
Moreover, no $1\neq g\in G$ fixes a curve   pointwise.
Thus it is enough to prove the following.

\begin{lem} Let $G$ be 
  an  Abelian group  acting
on a rational surface with cyclic stabilizers.
Assume also that no element fixes a curve of genus $\geq 1$ pointwise.
 Then
$G$ is  either $\z/m$,  $\z/m+\z/2$ or $\z/m+\z/3$.
\end{lem}

Proof.
We can take a $G$-equivariant resolution 
and then pass to the $G$-minimal model $T$.
(Note that if a smooth point $p\in T$ is fixed by
an Abelian group $H$ then $H$ also has a fixed point on
the blow up $B_pT$, so a fixed point on one
 model gives a fixed point on any other  model.)

Our aim is to find a $G$-invariant subset $Z\subset T$
such that either $Z\cong \p^1$ or $Z$ is at most 3 points.
In the first case, $G$ acts on $\p^1$ hence it  either
has a fixed point, or it acts through $(\z/2)^2$ and an index 2
subgroup has a fixed point. In the second case
a subgroup of index $\leq 3$ has a fixed point. Such
a subgroup is cyclic and so $G$ is
 $\z/m$,  $\z/m+\z/2$ or $\z/m+\z/3$ as required.

Consider first the case when  there is a 
$G$-equivariant ruling $f:T\to \p^1$. 
We are done if $f_*:G\to \aut(\p^1)$ is injective.
Otherwise, any $g\in \ker f_*$ fixes  2 points in each fiber of $f$,
thus $g$ fixes either 1 or 2 irreducible curves pointwise
and their union is $G$-invariant.
We are done if there is a $G$-invariant curve.

Otherwise, every $g\in \ker f_*$
 fixes 2 disjoint curves $E_1, E_2$ pointwise and $G$ permutes these curves.
By the Hodge index theorem 
$E_i^2\leq 0$, 
thus they generate the ``other'' extremal ray of the cone of curves.
In particular, either $T=\p^1\times \p^1$ 
or these curves are unique and  every element of
$\ker f_*$ fixes the same curves. If $G$ has a fixed point
 $p\in \p^1$ then $Z:=f^{-1}(p)\cap (E_1\cup E_2)$ is a 2--element set fixed
by $G$. If there is no such $p$ then  $G$ acts
on $\p^1$ and also on $E_1\cup E_2$ through $(\z/2)^2$.
Furthermore, there is an index 2 subgroup $H\subset G$ which
fixes each $E_i$ and has a fixed point on each $E_i$.
Thus $H$ is cyclic
 and so   $G\cong \z/m$ or   $G\cong \z/m+\z/2$.

The $\p^1\times \p^1$ case is left to the end.

Otherwise $T$ is a Del Pezzo surface of $G$-Picard number 1.

Assume that some $g\in G$  pointwise fixes some curves $C_i$.
 The $C_i$ are smooth,
 disjoint, rational curves.
By the adjunction formula, $C_i^2\in\{0,-1\}$ or $T\cong \p^2$.
Their sum $\sum C_i$ is $G$-invariant and not ample, but
this contradicts $G$-minimality. Thus either
$T\cong \p^2$ or every  element of $G$
acts with isolated fixed points.

Assume that  there is $H<G$ with  $H\cong (\z/2)^2$.
Then $T/H$ is a Del Pezzo surface with $A_1$-singularities only.
From Table 2 in (\ref{nontriv.pi1.lisat})
  we see that $\deg T/H=2$ and
$T\cong \p^1\times \p^1$. Thus, aside from this case, 
the 2--part of $G$ is cyclic.

 If $T=\p^2$,
take any $g\in G$. It has either $3$ fixed points or a fixed point 
and a fixed line (and the second case must happen if
$g^2=1$). In the second case, the fixed point is $G$-fixed
(and so $G$ is cyclic) 
and in the former case either all 3  points are $G$-fixed or
$G$ permutes them cyclically.

If $\deg T\in \{1,2,3,4,5\}$ then $\aut(T)$ acts faithfully
on $H^2(T,\z)$. Indeed, if $g$ acts trivially on $H^2(T,\z)$
then it descends to an automorphism of $\p^2$ with
$\geq 4$ fixed points in general position, hence $g=1$.
Thus, $G$ is a subgroup of the Weyl group of
$E_8,E_7,E_6,D_5,A_4$  (cf.\ \cite[Sec.25]{manin}).

If $\deg T=1$ then $|K_T|$ has a unique base point which is fixed by
$\aut(T)$, so $G$ is cyclic.

If $\deg T=2$ then there is a unique degree two morphism
$T\to \p^2$. If $|G|$ is odd and $H<G$ then 
 any $H$-fixed point on $\p^2$ is
dominated by (one or two) $H$-fixed point(s) on $T$. 
If $G$ is even, then a $G$-fixed point on $\p^2$
is  dominated by  two points and each is fixed by an
index 2 subgroup of $G$.

If $\deg T\in \{3,4,5\}$ then by looking at the order 
of the Weyl  groups, we see that the odd part 
$G_{odd}\subset \z/15$, except possibly
when $\deg T=3$ where the 3--part could be bigger.
Any action of a group on a cubic surface $T$ induces an
action on $\p^3\cong |K_T|$. For odd order groups this action lifts to
 $\c^4$ since the kernel of
$SL_4\to PSL_4$ is $\z/4$. Thus we get an eigenvector on $\c^4$
and a fixed point on the cubic $T$.

If $\deg T=6$ then only $\z/3$ can act on the $(-1)$-curves
nontrivially, and even for the $\z/3$-action there is
a $\z/3$-invariant set of 3 disjoint $(-1)$-curves. Thus the
action descends to $\p^2$.

There are no $G$-minimal  Del Pezzo surfaces of degree 7
and in degree 8 we get $\p^1\times \p^1$.

Thus we are left with $G$ acting on $\p^1\times \p^1$.
If any $g\in G$ interchanges the two factors, then
$g$ has 2 fixed points or a fixed rational curve
from the Lefschetz fixed point formula. In both cases 
an index 2 subgroup of $G$ has a fixed point.

The last case is when $G$ preserves the coordinate projections.
The image of $G$ in each $\aut(\p^1)$ is either cyclic
or $(\z/2)^2$. If the first case happens at least once, then
an index 2 subgroup has a fixed point. Finally we have to
deal with subgroups of $G\subset (\z/2)^4$. If the order is 8 or 16
then there are $g_1, g_2\in G$ which act trivially
on the first (resp. second) factor. Thus $\langle g_1,g_2\rangle$
is a noncyclic subgroup with a fixed point, a contradiction.
\qed
\medskip

Putting these together, we obtain the following.
Note that the torsion in $H_2$ is dual to the torsion in $H^3$,
thus a quotient of $H^3$ becomes a subgroup of $H_2$.

\begin{cor}\label{h2.tors.gen.thm} Let $L$ be a compact 5--manifold  
with a positive Sasakian structure.
Then  the torsion subgroup of the second homology, $\tors H_2(L,\z)$
has a subgroup $G$ such that
\begin{enumerate}
\item  $G$ is  $(\z/m)^2$,
 $(\z/5)^4$, $(\z/4)^4$,
  $(\z/3)^{2n}$ for $n\in \{2,3,4\}$ or $(\z/2)^{2n}$ and
\item $\tors H_2(L,\z)/G$ is 
 $\z/r$, $\z/r+\z/2$ or
 $\z/r+\z/3$ for some $r$. \qed
\end{enumerate}
\end{cor}

\begin{rem} Most Abelian groups can not be written in the above form.
If $G$ exists, then it is almost always unique up to
isomorphism. The only ambiguity is with the 6-torsion part.

To see this, we can consider the $p$-parts separately.
The main cases are
\begin{enumerate}
\item  $A_p=(\z/p^a)^2+\z/p^b$ and $A_p/G\cong \z/p^b$, 
\item  $A_p=\z/p^a+\z/p^b$ and $A_p/G\cong \z/p^{b-a}$.
\end{enumerate}
These are the only possibilities for $p\geq 7$,
with a few more cases for $p=2,3,5$.
\end{rem}

\begin{say}[Plan of the proofs of
(\ref{h2.tors.gen.cases}) and (\ref{main.thm})]\label{plan.say}
We follow the approach in \cite[6.8]{em5}. 
Start with  $(S,(1-\tfrac1{r})D+\Delta)$ satisfying the conditions
(\ref{algred.say}.1--3). Let $g:S'\to S$ be the
minimal resolution of $S$ and $h:S'\to S^m$ a minimal model of $S'$.
In the sequence of blow ups leading from $S$ to $S'$
and the subsequent blow downs leading from $S'$ to $S^m$
every intermediate surface $T$ satisfies the following condition:
\begin{enumerate}
\item[($*$)]  We can write $-K_T\equiv (1-\tfrac1{r})D_T+\Delta_T+H_T$
where  $D_T$ denotes the birational transform of $D$ on $T$,
 $\Delta_T$ is an effective linear combination
of rational curves (coming from $\Delta$ and the exceptional curves of $g$)
 and $H_T$ is nef and big
(this is a general divisor numerically equivalent to
 the pull back/push forward of 
$-(K_S+(1-\tfrac1{r})D+\Delta)$).
\end{enumerate}
This turns out to be very restrictive in many cases.

It is easy to see that $S$ is a rational surface, so
$S^m$ is either $\p^2, \p^1\times \p^1$ or a minimal ruled
surface $\f_n$ for some $n\geq 2$. For the latter, let
$E\subset \f_n$ denote the negative section and $F$ a fiber.
By an easy case analysis (cf.\ \cite[6.8]{em5}) we get the following.
\begin{enumerate}
\item If $r\geq 6$ and $g\geq 1$ then $S^m=\p^2, \p^1\times \p^1$
or $\f_2$ and $D^m\in |-K_{S^m}|$ is  smooth and elliptic.
\item If $r\geq 5, g\geq 2$ then $r=5$, $S^m=\f_3$ and $D^m\in |2E+6F|$
has genus $2$.
\item If $r\geq 3, g\geq 4$ then $r=3$, $S^m=\f_5$ and $D^m\in |2E+10F|$
has genus $4$.
\item If $r=4, g\geq 2$ then either $S^m=\f_3$ and $D^m\in |2E+6F|$
or $S^m=\f_2$ and $D^m\in |2E+5F|$. In both case $D^m$ has genus $2$.
\end{enumerate}
The plan is, in each case, to start with 
$(S^m,(1-\tfrac1{r})D^m+\Delta^m)$ and to write
down all possible $(S',(1-\tfrac1{r})D'+\Delta')$.
It is then easy to get a complete list of all
$(S,(1-\tfrac1{r})D+\Delta)$. Finally we need to check the additional
conditions, especially \cite[3.6 and 4.8]{em5}.
\end{say}

\section{The exceptional cases}

\begin{say}[The $(\z/5)^4$ and $(\z/3)^8$ cases]\label{5and3.pf}

Assume that $H_2(L.\z)$ contains a subgroup
isomorphic to $(\z/5)^4$ or $(\z/3)^8$.
Let $G\subset H_2(L.\z)$  be the subgroup given in
(\ref{h2.tors.gen.thm}). Since the odd torsion
in $H_2(L,\z)/G$ is either cyclic or  $(\z/3)^2$,
we conclude that either
\begin{enumerate}
\item $G= (\z/5)^4$, 
\item $G= (\z/3)^8$, or
\item  $G=(\z/3)^6$ and 
$H_2(L,\z)/G=(\z/3)^2$.
\end{enumerate}
Thus  $g(D)=2$ in the first case and $g(D)\geq 3$ in the second and third  
cases.

This implies that $S^0$ is simply connected.
Indeed, any nontrivial cover $\pi:S'\to S$ would give
another log Del Pezzo surface where $D':=\pi^{-1}(D)$
has genus $\geq 3$ in  the first case and genus $\geq 5$ in the other cases.
By the list in (\ref{plan.say}) there are no such surfaces.
In particular, the case  $G=(\z/3)^6$ and 
$H_2(L,\z)/G=(\z/3)^2$ does not happen.

Next, as in \cite[9.9]{em5}, we show that $S'=S^m$ and
$S$ is obtained by contracting the negative section $E$.
Indeed, any one point blow up of $\f_m$ maps
either to $\f_{m-1}$ (if we blow up a point not on $E$)
or  $\f_{m+1}$ (if we blow up a point  on $E$).
Since $S^m$ is unique in each case
by the list in (\ref{plan.say}),
$S'\to S^m$ is an isomorphism. 

If $S=S^m$ then 
$$
\Delta^m+H^m\equiv \sum (1-\tfrac1{r_i})D_i+H
$$
where the $D_i$ are rational and $H$ is ample. 
In the $(\z/5)^4$ case, this would give
$$
\tfrac25 E+\tfrac15F\equiv \sum (1-\tfrac1{r_i})D_i+H.
$$
The intersection number of the left hand side with
$F$ is $\tfrac25<\frac12$, so 
any $D_i$ is a fiber. Now intersecting with $E$ shows that
there  can not be any $D_i$. This is impossible since
 the left hand side is not  nef and big.

Similarly, in the $(\z/3)^8$ case we would need to solve
$$
\tfrac23 E+\tfrac13F\equiv \sum (1-\tfrac1{r_i})D_i+H,
$$
which is also impossible.

This proves  (\ref{main.thm}), except for the
equations. 

There is a systematic way to obtain $Y$ and the $\c^*$-action from
$(S,\Delta)$ (cf.\ \cite[2.5]{em5}), but I found it much easier to
do this by guessing. Note that by contracting the negative section
in $\f_n$, we get the weighted projective plane
$\p(1,1,n)$ and $D\in |2E+2nF|$ is a curve given by a
degree $2n$ equation. After completing the square,
these are of the form $z^2+f_{2n}(x,y)$.
The explicit relationship between $\Delta$ and the weights
exhibited in \cite[6]{bgk} now leads to 
$f=t^5+z^2+f_{6}(x,y)$ in case $r=5, S^m=\f_3$
and to 
$f=t^3+z^2+f_{10}(x,y)$ in case $r=3, S^m=\f_5$.\qed
\end{say}

The case analysis is harder for $(\z/4)^4$ and
we need some way to see how to get $S'$ from $S^m$.

\begin{say}[Blow up criterion]\label{bu.crit} Assume that
$\pi:T_1\to T$ is the inverse of the blowing up of
$p\in T$ with exceptional curve $E_1\subset T_1$ and that
$(T_1,(1-\tfrac1{r})D_1+\Delta_1+H_1)$ satisfies  (\ref{plan.say}.$*$).
That is, $-K_{T_1}\equiv (1-\tfrac1{m})D_1+\Delta_1+H_1$,
$\Delta_1$ is  effective  and $H_1$ is nef and big.

Set $D:=\pi_*D_1, \Delta:=\pi_*\Delta_1$ and
$H:=\pi_*H_1$. Then 
$(T,(1-\tfrac1{r})D+\Delta+H)$ also satisfies  (\ref{plan.say}.$*$).
Since $(K_{T_1}\cdot E_1)=-1$, we conclude that
$$
\left((1-\tfrac1{r})D_1+\Delta_1+H_1\right)\cdot E_1=1,
$$
hence 
$$
\mult_p \left((1-\tfrac1{r})D+\Delta+H\right)\geq 1.
$$
In practice we know $(T,(1-\tfrac1{r})D)$ 
 and we would like to
find $p$. This is possible if the numerical class of
$\Delta+H$ is small, but we get many possibilities if
$\Delta+H$ is bigger.
\end{say}

\begin{say}[The $(\z/4)^4$ case]
Here there are two possibilities for $(S^m, D^m)$
and $S'\to S^m$ need not be an isomorphism.

Let us start with  $S^m=\f_2$ and $D^m\in |2E+5F|$.
\medskip

{\it Case 1.}
If $S'=S_m$ then, as in (\ref{5and3.pf}),
 from $(E\cdot (\Delta^m+H^m))<0$ we see that
$E$ must be contracted. Thus $S$ is the 
 weighted projective plane
$\p(1,1,2)$ 
 and $D$ is given by a weighted degree 5 equation
$(f_5(x,y,z)=0)$ where $w(z)=2$.
(Alternatively,  $S$ is the quadric cone
in $\p^3$ and $D$ is a curve of degree 5
passing through its vertex.)
Thus we get the equation
$(t^4+f_5(x,y,z)=0)$ for $Y$.
This gives us (\ref{main.thm}.3.a).
\medskip

{\it Case 2.}
If we perform at least one blow up in going from $S^m$ to $S'$,
then as before, we could have ended up with
$S^m=\f_1$ or $S^m=\f_3$ instead.
The first case is impossible from the list of (\ref{plan.say})
 and the second one
we consider next.

\medskip
{\it Case 3.}
Thus assume that $S^m=\f_3$ and $C\in |2E+6F|$ is a
smooth curve. Thus  $\Delta^m+H^m\equiv \tfrac12 E+\tfrac12 F$.

If $S'=S^m$  then, as before,  we see that
$E$ must be contracted. 
This leads to the surface
$S=\p(1,1,3)$ and $D\in |\o_S(6)|$ a smooth curve.

This never leads to simply connected  5--manifolds  by \cite[4.8]{em5}.
All Seifert bundles are of the form
$Y(S, B, \frac{b}{4}D)$ where $B=aF$ for some $a\in\z$ and $b\in \{1,3\}$.
The corresponding Chern class is
$(a+\frac{3b}{2})F$, always a half integer.
Thus we get the basic cases with $c_1(Y/S)=\frac12$ and their
obvious quotients.

Since $D\in |\o_S(6)|$ and $6$ is even,
$H_1^{orb}(S,\frac34D)=\z/2$
and there are two basic cases: $a=-1,b=1$ and $a=-4,b=3$.

The double cover of $(S,\frac34D)$
is given by
$$
S_6=(f_6(x,y,z)+T^2=0)\subset \p(1,1,3,3),
$$
and the involution is $(x:y:z:T)\mapsto (x:y:z:-T)$.
Keep in mind that these are projective coordinates,
so the same  involution can be given as $(x:y:z:T)\mapsto (-x:-y:-z:T)$.
The basic Seifert bundle pulls back to
$Y=L^*(2,2,6,3;12)$ with typical equation
$$
x^6+y^6+z^2+t^4=0.
$$ 
There are 2 ways to lift the involution to a fixed point free
involution on $Y$. These are
$(x,y,z,t)\mapsto (-x,-y,-z,-t)$ and
$(x,y,z,t)\mapsto (-x,-y,-z,t)$.
These give the 2 families listed in  (\ref{main.thm}.3.b).

\medskip
{\it Case 4.}
Next we blow up at least 1 point on $S^m$.
This point can not be on $E$ since the resulting surface
would also dominate $\f_4$. 
Write
$$
\tfrac12 E+\tfrac12 F\equiv \Delta^m+H^m=\sum a_iF_i+bE+R^m,
$$
where the $F_i$ are distinct fibers and $R^m$
 has no irreducible component which is $E$ or a fiber.
By intersecting with $F$ we see that $b\leq \frac12$.
Intersecting with $E$ gives 
$$
0\leq (E\cdot R^m)= 3b-1-\sum a_i,
$$
in particular $b\geq \tfrac13$ and
$\tfrac12-b\leq \tfrac13(\tfrac12-\sum a_i)$,
thus $\sum a_i\leq \frac12$.
If $p$ lies on $F_1$ then
$$
\mult_p(\Delta^m+H^m)\leq a_1+(F\cdot R^m)
=a_1+(\tfrac12-b)\leq \tfrac16+\tfrac23 a_1-\tfrac13\sum_{i\neq 1}a_i
\leq \tfrac12.
$$
The condition (\ref{plan.say}.$*$) is 
$$
\mult_p(\tfrac34 D+\Delta^m+H^m)\geq 1,
$$
which is only possible if $p\in C$ and $a_1\geq 1/8$.
Furthermore, if $a_1\leq a_2$ then we get
$\tfrac 14\leq a_1\leq a_2$. But $\sum a_i\leq \tfrac12$,
and this would lead to $R^m=(\tfrac12-b)E$
which is impossible.
Thus we conclude:
\medskip

{\it Claim.} There is a unique fiber $F$ such that
$S'$ is obtained from $S^m$ by blowing up points
in or above $F\cap D$. 
\medskip

The general case is when $F\cap D$ consists of 2 points
and the special case is when $F\cap D$ is a single point
where $F$ and $D$ are tangent.

\medskip
{\it Case 5.}
Let us deal first with the general case and blow up $p\in F\cap D$.
We get $S_1\to S^m$ with exceptional curve $G_1$. Let
$F_1$ and $ D_1$ be the birational transforms of $F$ and $D$.
Write
$$
\Delta_1+H_1=aF_1+bE_1+cG_1+R_1.
$$
From (\ref{bu.crit}) we see that 
$$
c\leq \mult_p(\tfrac34 D+\Delta^m+H^m)-1
\leq \tfrac34+a+\tfrac13(\tfrac12-a)-1=\tfrac23a-\tfrac1{12}.
$$
So at every point  $q\in G_1\setminus(D_1+F_1)$ the multiplicity
of $\Delta_1+H_1$ is at most
$$
c+\mult_qR_1\leq c+\mult_pR^m\leq \tfrac23a-\tfrac1{12}
+ \tfrac16-\tfrac13 a=\tfrac1{12}+\tfrac13a<1,
$$
thus we can not blow these up. At $G_1\cap F_1$
the multiplicity
 is at most
$$
a+c+\mult_qR_1\leq a+\tfrac1{12}+\tfrac13a<1.
$$
Finally, at $G_1\cap D_1$
the multiplicity
 is at most
$$
\tfrac34+\tfrac1{12}+\tfrac13a=\tfrac56+\tfrac13 a<1,
$$
except when $a=\tfrac12$.
But  this again would lead to $R^m=(\tfrac12-b)E$
which is impossible.

Thus we conclude that we can blow up one or both of the points
$F\cap D$ and then we get $S'$. 

If we blow up only one point, we have to contract $F_1$
and we are in the already considered case when $S'=\f_2$.

\medskip
{\it Case 6.}
If we blow up both points of $F\cap D$, we have to contract $E_1\cup F_1$.
We get a surface with Picard number 2 and a single
cyclic quotient singularity of the form
$\c^2/\tfrac15(1,2)$.

\medskip

{\it Claim.} The surface $S$ is isomorphic to a
quasi-smooth hypersurface $S_6\subset \p(1,1,3,5)$
and $D$ is the complete intersection of $S_6$ with 
$(t=0)$.
\medskip

I found this isomorphism by computing the quotient 
by the hyperelliptic involution of $D$ which acts
on $S^m,S'$ and also on $S$. Once the isomorphism is guessed,
it is easier to verify  by working backwards from 
the surface $S_6\subset \p(1,1,3,5)$.
Its equation can be written, after coordinate changes, as
$$
S_6=(z^2+f_6(x,y)+\ell_1(x,y)t=0)\subset \p(1,1,3,5).
$$
Notice that its intersection with $(\ell_1=0)$
is the reducible curve $(z^2+f_6(x,\alpha x)=0)\subset \p(1,3,5)$
for some $\alpha$.
The two irreducible components correspond to the two 
exceptional curves of $S'\to S^m$.
The rest is a straightforward computation. 
The surface is specified by  choosing 
6 points in $\p^1$ (given by $f_6=0$) plus one more
corresponding to the choice of $\ell_1$.

This leads to the case (\ref{main.thm}.3.e).

\medskip
{\it Case 7.}
Next we deal with the special case  when $F\cap D$ is a single point
where $F$ and $D$ are tangent. Computations as above
yield that in this case we can blow up the point $p$ on $D$ 
at most 3--times.

\medskip
{\it Case 8.}
If $S'$ is obtained by 1 blow up, we factor through $\f_2$ as before.
If we do 2 blow ups, we get $S=\p(1,2,5)$ and
$D\in |\o(10)|$.  
As in Case 3, $H_1^{orb}(S,\frac34D)=\z/2$.
Set $L:=(x=0)$.
There are two basic cases, corresponding to
$Y(S,-2L, \frac14D)$ and $Y(S,-7L, \frac34D)$.
The index 2 point $(0:1:0)\in S$ adds a further complication
since the first of these does not satisfy  the smoothness  condition
\cite[4.8.1]{em5}.
 
The orbifold double cover of $(S,\frac34D)$ is
$$
S_{10}:=(f_{10}(x,y,z)+T^2=0)\subset \p(1,2,5,5)
$$
and the involution is
$(x:y:z:T)\mapsto (x:y:z:-T)$.
Since these are weighted projective coordinates, this is the same as
$(x:y:z:T)\mapsto (-x:y:-z:T)$.
(Note that $-1$ acts by sending a coordinate $u$ to $(-1)^{wt(u)}u$,
hence the $+$ sign in front of $y$.)

This leads to the basic examples
$L=L^*(2,4,10,5;20)$ with sample equation  $(x^{10}+y^{5}+z^2+t^4=0)$,
The lifting of the involution is
either $(x,y,z,t)\mapsto (-x,y,-z,-t)$ or
 $(x,y,z,t)\mapsto (-x,y,-z,t)$. Here the second action has fixed point,
this is consistent with our earlier considerations.
Thus we get only one family, as in  (\ref{main.thm}.3.d).

\medskip
{\it Case 9.}
Finally, if we blow up 3--times, we can contract the
birational transforms of $E,F$ and of the first 2
exceptional curves. This gives a  surface with a single singular point
of the form $\c^2/\tfrac19(2,5)$.

\medskip

{\it Claim.} The surface $S$ is isomorphic to a
quasi-smooth hypersurface $S_{10}\subset \p(1,2,5,9)$
and $D$ is the complete intersection of $S_{10}$ with 
$(t=0)$.
\medskip

 Once again, 
it is easier to verify this by working backwards.
The equation of $S_{10}$ is
$$
(f_5(x^2,y)+z^2+xt=0)\subset \p(1,2,5,9).
$$
Its intersection 
with $(x=0)$ is the  curve $(ay^5+bz^2=0)\subset \p(2,5,9)$.
This is a smooth (but not quasi-sooth) rational curve,
and it  corresponds to the last 
exceptional curve of $S'\to S^m$.
The surface is specified by  choosing 
6 points in $\p^1$ (given by $xf_5=0$), that is
a genus 3 hyperelliptic curve plus a specified branch point.
The rest is a straightforward computation. 

Thus we obtain (\ref{main.thm}.3.c).
\qed
\end{say}

\begin{say}[Other examples] It is easy to write down
infinitely many positive Sasakian structures
on  certain simply connected 5--manifolds $L$
either by hand or by consulting the
(partially unpublished) lists of Boyer and Galicki.
Below I write the orbifolds $(S,\Delta)$ and the simplest equation.
(Here $C_k$ denotes a general curve of degree $k$ 
in the weighted projective plane $\p(a,b,c)$ with coordinates
$x,y,z$ and
$\ell:=(x=0)$.)
\begin{enumerate}
\item $H_2(L,\z)=(\z/3)^2$. For any $(k,6)=1$, take 
$$
\bigl(\p(1,1,2), (1-\tfrac13)C_4+(1-\tfrac1{k})\ell\bigr)\qtq{and}
x^{4k}+y^4+z^2+t^3.
$$
\item $H_2(L,\z)=(\z/5)^2$. For any $(k,30)=1$, take 
$$
\bigl(\p(1,2,3), (1-\tfrac15)C_6+(1-\tfrac1{k})\ell\bigr)\qtq{and}
x^{6k}+y^3+z^2+t^5.
$$
\item $H_2(L,\z)=(\z/3)^4$. For any $(k,30)=1$, take 
$$
\bigl(\p(1,2,5), (1-\tfrac13)C_{10}+(1-\tfrac1{k})\ell\bigr)\qtq{and}
x^{10k}+y^5+z^2+t^3.
$$
\item $H_2(L,\z)=(\z/2)^{2n}$. For $n\geq 0$ and $(k,2n(2n+1))=1$  take
$$
\bigl(\p(1,1,n), (1-\tfrac12)C_{2n+1}+(1-\tfrac1{k})\ell\bigr)\qtq{and}
x^{k(2n+1)}+y^{2n+1}+yz^2+t^2.
$$
\end{enumerate}
Sasakian--Einstein structures are known to exist
in the first 3 cases, but for the last one the  criterion
of \cite{dk} fails and existence is not known.

The situation is more complicated in the next two cases:
\begin{enumerate}\setcounter{enumi}{4}
\item $H_2(L,\z)=(\z/4)^2$ or $(\z/6)^2$. For any $k$, take 
$$
\bigl(\p(1,2,3), (1-\tfrac14)C_6+(1-\tfrac1{k})\ell\bigr)
\qtq{resp.}
\bigl(\p(1,2,3), (1-\tfrac16)C_6+(1-\tfrac1{k})\ell\bigr).
$$
These are the right examples as log Del Pezzo surfaces, but
the condition \cite[4.8.2]{em5} fails, so the corresponding Seifert bundles
are not simply connected. Since my approach is to rule out everything at the
surface level, these cases could be very hard to settle using my methods.
\end{enumerate}

 This leaves
open the finiteness question for
\begin{enumerate}\setcounter{enumi}{5}
\item $H_2(L,\z)=(\z/m)^2$ with  $7\leq m\leq 11$, and
\item $H_2(L,\z)=(\z/3)^6$.
\end{enumerate}
The  computations in Section \ref{klt.sect} suggest that in these cases
there should be only finitely many families of 
pre-SE  structures.

However, as the examples with $H_2(L,\z)=(\z/4)^4$  suggest, the
case analysis can be rather tricky and unexpected special
configurations may arise.
\end{say}

\section{The main series}

Let us start with fixing an error in \cite{em5}.

\begin{say}[Correction to {\cite[Thm.1.8.1(b)]{em5}} ]
\label{corr181.say}
The theorem asserts that for each $m\geq 12$, there are
exactly 93  pre-SE Seifert bundles
$f:L\to (S,\Delta)$ on 
 compact 5--manifolds $L$  
satisfying $H_1(L,\z)=0$ and $\tors H_2(L,\z)\cong (\z/m)^2$. 

The construction of the 93 families of pairs
$(S,\Delta)$ is correct.
Going from $(S,\Delta)$ to the Seifert bundle is, however,
done incorrectly  since the conditions for  a
Sasakian structure and for  a pre--Sasakian--Einstein
structure have been thoroughly mixed up.

The construction of the 93 families given in \cite[7.6]{em5}
starts with a surface $T$ which is
one of $\p^1\times \p^1, \p^2, Q,  S_5, \p(1,2,3)$.
For each of these write $-K_T\sim d(T)H$ where $H\in \weil(T)$
is a positive generator. We have
$$
d(\p^1\times \p^1)=2,  d(\p^2)=3, d(Q)=4,  d(S_5)=5,  d(\p(1,2,3))=6.
$$
Next we perform some weighted blow ups \cite[7.3]{em5}
to get $S=B_{m_1,\dots,m_k}T$.
There are $k$ exceptional curves $E_1,\dots, E_k$.
Each $E_i$ passes through a unique singular point $p_i\in S$
and $E_i$ generates the local class group which is $\z/m_i$.
Set $d(S):=\gcd(m_1,\dots,m_k, d(T))$.

The divisor class group $\weil(S)$ is freely generated by
$$
\pi^*H, E_1,\dots, E_k\qtq{and}
K_S\sim -d(T)\pi^*H+\sum m_iE_i.
$$

The main condition that was overlooked is the
smoothness criterion for Seifert bundles \cite[3.6]{em5}.

In our case there is only one curve $D$ and $S$ is smooth along
$D$. Thus, by \cite[3.6]{em5},
 the corresponding Seifert bundle 
$Y(S,B,\frac{b}{m}D)$ is smooth iff 
$$
\mbox{$B$ generates the local class group at every  point.}
\eqno{(\ref{corr181.say}.1)}
$$
\medskip

{\it Claim \ref{corr181.say}.2.} Notation as above. Then
$$
Y(B_{m_1,\dots,m_k}T, aH+c_1E_1+\cdots+c_kE_k, \frac{b}{m}D)
\qtq{is smooth}
$$
iff  $(m_i,c_i)=1$ for $i=1,\dots,k$ and
$(a,d(T))=1$ if $T$ is singular.
\medskip

As a consequence, we see that all 93 cases correspond 
to positive Sasakian structures on smooth 5--manifolds.

A Seifert bundle is pre-SE iff
its Chern class $c_1(Y/S)=B+\frac{b}{m}D$ is a rational multiple of 
$-(K_S+(1-\frac1{m})D)$. In our case $D\sim -K_S$, hence
$B$ itself is a rational multiple of $-K_S$.
Thus
$$
B=r\bigl(\tfrac{d(T)}{d(S)}\pi^*H-\sum \tfrac{m_i}{d(S)}E_i\bigr)
$$
for some positive integer $r$.
Thus if  (\ref{corr181.say}.1)
holds then $d(S)=m_i$ for every $i$.

Furthermore, in the cases when $T$ is singular,
$B$ generates the local class group at each singular point of $T$
iff $d(T)=d(S)$. 
Thus we obtain the following.
\medskip

{\it Claim \ref{corr181.say}.3.} Notation as above. Then
$$
Y(B_{m_1,\dots,m_k}T, aH+c_1E_1+\cdots+c_kE_k, \tfrac{b}{m}D)
\qtq{is smooth and pre--SE iff}
$$
\begin{enumerate} 
\item $m_1=\cdots=m_k=d(S)$, 
\item 
$aH+c_1E_1+\cdots+c_kE_k= r(-\frac{d(T)}{d(S)}H+E_1+\cdots+E_k)$
for some $r\in \z$, and
\item  if $T$ is singular, then $(r,d(T))=1$ and $d(S)=d(T)$. \qed
\end{enumerate}
\medskip

These conditions cut down considerably
the list given in \cite[7.6]{em5} and we
 get the following 19  cases:
$$
\begin{array}{c}
\begin{array}{lc}
\qquad \mbox{surfaces} & \mbox{$d(S)$}\\*[1ex]
B_{1}\p^2, B_{11}\p^2,\dots, B_{11111111}\p^2 & 1\\
\p^1\times \p^1, 
B_2\p^1\times \p^1, B_{22}\p^1\times \p^1, B_{222}\p^1\times \p^1 & 2\\
\p^2, B_3\p^2, B_{33}\p^2 &  3\\
Q, B_4Q & 4\\
S_5 & 5\\
\p(1,2,3) & 6
\end{array}
\ \\
\mbox{Table 1}
\end{array}
$$

The condition \cite[4.8.2]{em5} says that
the resulting Seifert bundle $L$ satisfies
$H_1(L,\z)=0$ iff 
$H^2(S,\z)\to H^2(D,\z)\to \z/m $ is surjective.
If $S=T$, that is, we do no blow ups at all, then 
by \cite[9.8]{em5}
$$
\im\bigl[H^2(T,\z)\to H^2(D,\z)\bigr]=d(T)H^2(D,\z).
$$
After blow ups, the new curves that come in are the
 exceptional curves $E_i$. Here $(E_i\cdot D)=1$ but the $E_i$ pass through the
singular points and they are only homology classes.
To get a cohomology class (or Cartier divisor), we need to
take $m_iE_i$ since $m_i$ is also the index of the singular point.
Thus we conclude:
\medskip

{\it Claim \ref{corr181.say}.4.} For $S=B_{m_1,\dots,m_k}T$,
$\im\bigl[H^2(S,\z)\to H^2(D,\z)\bigr]= d(S)\cdot
H^2(D,\z)$.\qed
\end{say}

\medskip

\begin{cor} Let $S$ be a projective surface with 
Du Val singularities 
such that $\pi_1(S^0)=1$. Let
 $D\in |-K_S|$ be a smooth elliptic curve.
There is a simply connected Seifert bundle
$f:L\to (S,(1-\frac1{m})D)$ with a pre-SE structure iff
\begin{enumerate}
\item $S$ is one of the surfaces in Table 1, and
\item  $m$ is relatively prime to $d(S)$.
\end{enumerate}
In these cases, $f:L\to (S,(1-\frac1{m})D)$ is
uniquely determined by $(S,(1-\frac1{m})D)$
(up to reversing the orientation of the fibers) and
$L$ carries a Sasakian--Einstein metric for $m\geq 7$.
\end{cor} 

Proof. we have proved eveything, except the claims
about the existence of Sasakian--Einstein metrics.
This will be established in (\ref{ke.conds}).\qed

\begin{say}[Equations] Quite surprisingly, all the singular surfaces on the
list can be realized in weighted projective 3--spaces.
All of these examples are on the Boyer--Galicki lists.
Here  also claim the converse: every singular Del Pezzo surface
in Table 1 is isomorphic to a corresponding surface below.
\begin{enumerate}
\item $B_2\p^1\times \p^1$:  $S_3\subset \p(1,1,1,2)$ with simplest equation
$$
x^3+y^3+z^3+xt^m=0 \qtq{for $(m,2)=1$.}
$$
\item $B_{22}\p^1\times \p^1$:  $S_4\subset \p(1,1,2,2)$ with simplest equation
$$
x^4+y^4+z^2+zt^m=0 \qtq{for $(m,2)=1$.}
$$
\item $B_{222}\p^1\times \p^1$:  $S_6\subset \p(1,2,3,2)$
 with simplest equation
$$
x^6+y^3+z^2+t^{3m}=0 \qtq{for $(m,2)=1$.}
$$
\item $B_3\p^2$:   $S_4\subset \p(1,1,2,3)$ with simplest equation
$$
x^4+y^4+z^2+xt^m=0 \qtq{for $(m,3)=1$.}
$$
\item $B_{33}\p^2$:  $S_6\subset \p(1,2,3,3)$ with simplest equation
$$
x^6+y^3+z^2+zt^m=0 \qtq{for $(m,3)=1$.}
$$
\item $Q$:  $S_4\subset \p(1,1,2,4)$ with simplest equation
$$
x^4+y^4+z^2t^m=0 \qtq{for $(m,2)=1$.}
$$
\item $B_4Q$:  $S_6\subset \p(1,2,3,4)$ with simplest equation
$$
x^6+y^3+z^2+yt^m=0 \qtq{for $(m,2)=1$.}
$$
\item $S_5$:  (sorry for the notation) 
$S_6\subset \p(1,2,3,5)$ with simplest equation
$$
x^6+y^3+z^2+zt^m=0 \qtq{for $(m,5)=1$.}
$$
\item $\p(1,2,3)$:  $S_6\subset \p(1,2,3,6)$ with simplest equation
$$
x^6+y^3+z^2+t^m=0 \qtq{for $(m,6)=1$.}
$$
\end{enumerate}
\end{say}

\begin{rem} There are no hypersurface links 
on the Boyer--Galicki lists giving infinitely many
$(S,\Delta)$ where $S\in \{B_1\p^2, \dots, B_{11111}\p^2\}$.

I claim that these  can not be realized
 as the link $L$ of a hypersurface with a $\c^*$-action, 
at least when $H_2(L,\z)\supset (\z/m)^2$ for
$m\geq 12$ and $(S,\Delta)$ is log Del Pezzo.
Assume the contrary. Then we get that
$S$ is a hypersurface in a weighted projective space
$\p(a,b,c,d)$ such that $K_S+(1-\frac1{m})D$ is proportional to
$H|_S$ where $H$ is the hyperplane class 
of the weighted projective space. Since $D\in |-K_S|$, we conclude that
$K_S$ is proportional to
$H$.  By the Grothendieck--Lefschetz theorem,
this implies that $-K_S=dH$ for some $d\in \z$. 
In our cases $-K_S$ is not divisible, so $d=1$.
This implies that
$$
h^0(S,\o_S(-K_S))=h^0(\p(a,b,c,d), \o_{\p}(1))\leq 4,
$$
since this dimension is the number of times that 1 occurs among $a,b,c,d$.
In the above cases, however, $h^0(S,\o_S(-K_S))\geq 5$.

Galicki told me that the link $L^*(2,3,4,7;14)$
realizes  $S=B_{1111}\p^2$, but the corresponding $\Delta$
is not a rational multiple of $-K_S$.
\end{rem}

\begin{say}[The cases with nontrivial fundamental group]
\label{nontriv.pi1.lisat} The classification of minimal 
Del Pezzo surfaces with Du Val singularities
is completed in \cite{furushima, mi-zh1, mi-zh2, ye}. 
We are interested only in those that have cyclic quotient singularities.
The 5 cases where $\pi_1(S^0)=1$ were considered in \cite{em5}.
The following table lists the remaining ones.
$$
\begin{array}{c}
\begin{array}{cccccc}
\mbox{degree} & \mbox{singularities}& \pi_1(S^0) & \weil/\pic & 
\mbox{univ. cover}& \pi_1(S^0\setminus D)\\*[1ex]
1 & A_8 & \z/3 & \z/3 & B_{3111}\p^2 & \z/3\\
1 & A_7+A_1  & \z/4  &\z/4 & B_{2111}\p^2& \z/4\\
1 & A_5+A_2+A_1 & \z/6 & \z/6 &B_{111}\p^2& \z/6\\
1 & 4A_2 & (\z/3)^2 & (\z/3)^2 & \p^2& G_{27}\\
1 & 2A_3+2A_1 & \z/2+\z/4 & \z/2+\z/4 &\p^1\times \p^1 & G_{16}\\
1 & 2A_4 & \z/5 & \z/5 &B_{1111}\p^2 & \z/5\\
1 & A_3+4A_1 & (\z/2)^2 & (\z/2)^3+\z/4 & B_{22}\p^1\times \p^1  & \z/2+\z/4\\
2 & A_7  & \z/2 & \z/4 & B_{41}\p^2& \z/2\\
2 & A_5+A_2 & \z/3 & \z/6 &B_{11}Q & \z/3\\
2 & 2A_3+A_1 & \z/4 & \z/2+\z/4 &\p^1\times \p^1& \z/2+\z/4 \\
2 & 6A_1 & (\z/2)^2 & (\z/2)^4 &\p^1\times \p^1& G_8 \\
2 & 2A_3 & \z/2 & \z/2+\z/4 & B_{22}\p^1\times \p^1 & \z/4\\
3 & A_5+A_1 & \z/2 & \z/6 & B_3\p^2 & \z/6\\
3 & 3A_2 & \z/3 & (\z/3)^2 & \p^2& (\z/3)^2\\
4 & A_3+2A_1 & \z/2 & \z/2+\z/4 & Q & (\z/2)^2 \\
4 & 4A_1 & \z/2 & (\z/2)^3 & \p^1\times \p^1 & (\z/2)^2\\
\end{array}\\
\ \\
\mbox{Table 2}
\end{array}
$$
Here $G_n$ is a nonabelian group of order $n$.
\begin{enumerate}
\item $G_8$ is the quaternion group,
\item $G_{16}\subset GL_4$ is generated by
$(x,y,z,t)\mapsto (z,t,-x,y)$ and 
$(x,y,z,t)\mapsto (y,-x,t,-z)$.
\item $G_{27}\subset GL_3$ is generated by
$(x,y,z)\mapsto (y,z,x)$ and 
$(x,y,z)\mapsto (x,\epsilon y,\epsilon ^2 z)$ where $\epsilon^3=1$.
\end{enumerate}

The computation of the table:
The papers \cite[p.13-15]{furushima}, \cite[p.71]{mi-zh1}, 
\cite[p.193]{mi-zh2} and \cite{ye}
contain tables for the first 4 columns, except
$\weil/\pic$ in the 4 cases with Picard number 2.  (The fundamental group of
$A_7+A_1$ is listed in \cite[p.71]{mi-zh1} as $(\z/2)^2$, but it
is $\z/4$.) The Picard number and the singularities of the
universal cover are listed in \cite[p.71]{mi-zh1}, from this it
is easy to work out where the surface is on the list \cite[7.6]{em5}.

By \cite[1.2 and 1.6]{ye}, for each singularity type there is a unique surface,
except for $2A_3$ for which there is a 1--parameter family.

Once we have $\bar S\to S$ as the universal cover
and $\bar D\subset \bar S$ is the preimage of $D$ then we have
an exact sequence
$$
\pi_1(\bar S^0\setminus \bar D)\to 
\pi_1(S^0\setminus  D)\to \pi_1(S^0)\to 1.
$$
Each time $\bar S$ is obtained by a blow up of weight 1,
the resulting $\p^1\subset \bar S$ intersects $\bar D$ transversally
at 1 point, so $\pi_1(\bar S^0\setminus \bar D)=1$.
In the remaining cases one needs to write down the
action of the $\pi_1(S^0)$ and see how it lifts to
the universal cover of $\bar S^0\setminus \bar D$.

Computing any entry of the table is an elementary task.
Some computations are quick but a few are quite tedious.
It is unfortunately easy to miss or misdraw a $-1$-curve
after performing many blow ups, so anyone wishing to rely on
a particular entry is advised to recheck it.
\end{say}

In  the simply connected case there are many isomorphisms between
blow ups, but this does not happen for the general case.

\begin{lem}\label{no.two.isom}  Let $S$ be a Del Pezzo surface
 with cyclic quotient singularities such that $|\pi_1(S^0)|> 1$.
There is a unique line in Table 2 such that $S$ is a weighted blow up
of a surface on that line.
(We do not claim that the blow up itself is unique.)
\end{lem}

Proof. The blow ups do not change the fundamental group
\cite[7.3]{em5} and we create only $A_1$ and $A_2$ singularities
since we can only have blow ups $B_{m_1,\dots,m_r}T$
where $\sum m_i<\deg T-1\leq 3$.
It turns out that the fundamental group and the collection
of $A_i:i\geq 3$ singularities uniquely determine in which line
of  Table 2 the surface is. The Picard number and the
$A_i:i\leq 2$ singularities now determine the number and
type of  blow ups performed.\qed
\medskip

If $\deg T=1$ (resp.\ $2,3,4$) then we get $1$ (resp.\ $2,4,7$) blown up 
surfaces,
including $T$ itself. 
Thus we get
39 deformation types of  Del Pezzo surfaces
 with cyclic quotient singularities such that $|\pi_1(S^0)|> 1$.
The 93 cases where $\pi_1(S^0)=1$ were enumerated
in \cite[7.6]{em5}, giving
 a total of 132 deformation types.

\begin{say}[Existence of smooth Seifert bundles] \label{sm.sb.ex.say}
Let $S$ be one of the surfaces in Table 2
and $D\in |-K_S|$ a smooth elliptic curve.
As in (\ref{corr181.say}),
we are considering Seifert bundles
$Y(S,B,\frac{b}{m}D)$ where $B$ is a Weil divisor class on $S$.

By (\ref{corr181.say}.1) $Y(S,B,\frac{b}{m}D)$ is smooth iff 
 $B$ generates the local class group at every  point.

Finding all such $B$ 
 requires a detailed computation of the divisor class group
$\weil(S)$ and the restriction map
$$
\weil(S)/\pic(S)\to \sum_{p\in \sing S} \weil(p,S).
$$
On a Del Pezzo surface of degree $\leq 7$, the
curves $C$ with $(C\cdot K_S)=-1$ generate $\weil(S)$.
On the minimal desingularization $S'\to S$ these are the $-1$ curves.
Thus if we have a description of $S'$ as a blow up
of $\p^2$, we see all such curves by looking at
lines through 2 blow up points, conics through 5 blow up points,
etc. (See \cite[Sec.26]{manin} for the complete list in degrees 2 and 1.)
The description given in \cite{furushima} gives exactly these
blow ups. See also \cite[Figure 1]{mi-zh1}.

In some cases 
$\weil(S)/\pic(S)$ is too small to get surjection onto
some $\weil(p,S)$, and then  
 there are no smooth Seifert bundles at all.
This happens in 3 cases:
$$
A_8, A_7+A_1, A_7.
$$
More surprising is the mildly singular $A_3+2A_1$ case
which again has no smooth Seifert bundle over it.
On the minimal desingularization, the configuration of
$-1$ and $-2$ curves is
$$
\stackrel{-2}{\circ} - \stackrel{-1}{\circ} - \stackrel{-2}{\circ} - 
\stackrel{-2}{\circ} - \stackrel{-2}{\circ} - \stackrel{-1}{\circ} - 
\stackrel{-2}{\circ}
$$
In order to generate both $\z/2$ on the ends, we need to take the
two $-1$-curves with odd coefficients, but then
we get only twice the generator of the $\z/4$
of the middle singularity.

The more complicated $A_5+A_1$ case leads to the configuration
$$
\stackrel{-2}{\circ} - \stackrel{-1}{\circ} - \stackrel{-2}{\circ} - 
\stackrel{-2}{\circ} - \stackrel{-2}{\circ} - 
\stackrel{-2}{\circ} - 
\stackrel{-2}{\circ}
$$
plus an extra $-1$-curve. The $-1$-curve shown
generates both local class groups, so there are
smooth Seifert bundles, even with SE structure. 

A glance at the diagrams for $A_5+A_2+A_1$ and for
$2A_4$ in \cite[Figure 1]{mi-zh1} shows $-1$-curves
which generate all local class groups.

As another concrete example, 
the group $G_{27}$ operates freely on 
$\c^3$ outside the origin, thus
$S^5/G_{27}\to \p^2/(\z/3)^2$ is a smooth Seifert bundle.
\end{say}

Just for illustration, let us compute one simple case
completely.

\begin{exmp}[$3A_2$ case]  We can write this as
$S=\p^2/(\z/3)$ by the action 
$$
(x:y:z)\mapsto (x:\epsilon y: \epsilon^2 z)\qtq{where} \epsilon^3=1.
$$
We can take $D=(x^3+y^3+z^3=0)$. The universal cover of
$\p^2\setminus D$ is the cubic
$$
(x^3+y^3+z^3=t^3)\subset \p^3.
$$
We get a $(\z/3)^2$-action generated by
$$
(x:y:z:t)\mapsto (x:\epsilon y: \epsilon^2 z:t)\qtq{and} 
(x:y:z:t)\mapsto (x: y: z:\epsilon t).
$$

The 3 coordinate lines give the curves $A,B,C\subset \p^2/(\z/3)$.
These generate $\weil(S)$ subject to the relations
$3A=3B=3C=A+B+C$. We can thus rewrite
$$
\weil(S)=\z[A]+\z/3[A-B].
$$
By explicit computations, the class
$B_{uv}:=uA+v(A-B)$ corresponds to a smooth Seifert bundle
iff
none of the numbers $u, v, u+v$ is divisible by 3.
Whenever this holds, there is a smooth pre-SE  Seifert bundle
$L_{uv}\to (S, B_{uv}, (1-\frac1{m})D)$ for any $(m,3)=1$.
It has a Sasakian--Einstein metric for every $m\geq 4$.
\end{exmp}

\begin{rem} As a consequence of the classification,
\cite[Thm, p.184]{mi-zh2} concludes that the fundamental group
of a  Del Pezzo surface with Du Val singularities is abelian.

This is easy to see directly as follows. Let $g:S'\to S$ be
the universal cover. Pick a smooth elliptic curve
$C\in |-K_S|$. Then $C':=g^{-1}(C)\in |-K_{S'}|$, thus
it is also a smooth elliptic curve. 
Hence the fundamental group is the same as the kernel
of the group homomorphism $C'\to C$, hence an abelian group
with at most 2 generators.
\end{rem}

\section{Klt conditions}
\label{klt.sect}

While we did not use it in the proof, it is instructive to see
the relationship between finding $S'$ and the
klt condition for $(S^m,(1-\tfrac1{r})D^m+\Delta^m+H^m)$.

Since $\Delta'$ is a nonnegative linear combination, all
exceptional curves of $S'\to S^m$ have nonpositive
discrepancy with respect to
$(S^m,(1-\tfrac1{r})D^m+\Delta^m+H^m)$. If the latter
pair is klt, then, by an observation of Shokurov,
there are only finitely many such curves on
any birational model of $S^m$ and they can be found
explicitly (cf.\ \cite[2.36]{kmbook}). 

The problem with using this result is that
 we do not know $\Delta^m$ and $H^m$, only
the numerical class 
$$
 \Delta^m+H^m\equiv -\left(K_{S^m}+(1-\tfrac1{r})D^m\right).
$$
The usual proofs of the above finiteness result start
by taking a log resolution of 
$(S^m,(1-\tfrac1{r})D^m+\Delta^m+H^m)$, and in essence
we would need to understand all smooth surfaces
dominating $S^m$.

The key part of the proofs 
in (\ref{5and3.pf}) and (\ref {corr181.say}) is to describe all
exceptional curves over $S^m$ which have nonpositive
discrepancy with respect to some
$(S^m,(1-\tfrac1{r})D^m+\Delta^m+H^m)$.

To this end one can use 
the following result which sharpens the klt conditions
used in \cite[Sec.8]{em5}.

\begin{prop} \label{klt.crit} Let $S$ be a smooth surface,
$C$ a smooth curve on $S$ and $D$ an effective $\q$-divisor
on $S$ such that $C\not\subset \supp D$.
 Let $p\in S$ be  point and $n$ a natural number.
Then $(S, (1-\frac1{n})C+cD)$ is klt at $p$ for
$$
c< \min\left\{
\frac1{(C\cdot D)_p}+\frac1{n\cdot \mult_pD},\frac1{\mult_pD}\right\}.
$$
\end{prop}

Proof. Choose local coordinates $(x,y)$ such that
$C=(y=0)$ and let $f(x,y)=0$ be an equation of $mD$ for some $m$
such that $mD$ is an integral divisor.

Consider the local degree $n$ cover $\pi:T\to S$ given by $y=z^n$.
By \cite[20.3.2]{k-etal},  $(S, (1-\frac1{n})C+cD)$ is klt at $p$
iff $(T, \frac{c}{m}(f(x,z^n)=0))$ is klt at $p$.

We aim to apply a theorem of Var\v cenko \cite{varc}
which gives a condition for  $(T, \gamma(g(x,z)=0))$ to be klt
in terms of the Newton polygon of $g$ in a suitable coordinate system
which is achieved after a series of coordinate changes
of the form $(x,z)\mapsto (x-\alpha z^i, z)$ or
$(x,z)\mapsto (x, z-\alpha x^i)$.
The problem is that we can handle only those coordinate changes
which are compatible with $\pi$. 
That is, those of the form $(x,z)\mapsto (x-\alpha z^{ni}, z)$. 
Thus we have to look
carefully at the proof, not just the final result.

We state the result in the rather artificial form that we 
need.
 The reader should consult the proof given in
\cite[6.40]{ksc}, especially pp.172--3.

\begin{lem} \label{varch.thm}
Write $g(x,z)=\sum b(i,j)x^iz^j$. Assume that one of the
following holds.
\begin{enumerate}
\item (Main case.) There are $(i,j)$ and $(i',j')$ such that $b(i,j)\neq 0$,
$b(i',j')\neq 0$ and the line segment
$[(i,j), (i',j')]$ contains a point $(\gamma,\gamma)$ with $\gamma< c^{-1}$.
\item (Degenerate case.) There is $(i,j)$  such that $b(i,j)\neq 0$
and $i,j<c^{-1}$.
\end{enumerate}
Then one of the following also holds:
\begin{enumerate}\setcounter{enumi}{2}
\item (Klt case.) $(T, c(g(x,z)=0))$ is klt.
\item (Coordinate change.) There are natural numbers
$u,v,w$ and $e$  such that $(u,v)=1$,  $e>w/(u+v)$ and 
$$
\sum_{ui+vj=w} b(i,j)x^iy^j\qtq{is divisible by}
(\alpha x^v+\beta y^u)^e.
$$
In this case necessarily $u=1$ or $v=1$, there is a unique such factor
and the  new coordinates are
$(x,\alpha x^v+\beta y)$ if $u=1$ and
$(\alpha x+\beta y^u,y)$ if $v=1$.
\end{enumerate}
\end{lem}

Set $a=(L\cdot D)_p$ and $d=\mult_pD$.
As before,  write
$f=\sum b(i,j)x^iy^j$.
Then $b(am,0)\neq 0$ and
$b(i,j)\neq 0$ for some $i+j=md$.
Thus
$$
f(x,z^n)= \sum b(i,j)x^iz^{nj}.
$$
If $b(i,j)\neq 0$ for some $i+j=md$ and $j\geq md/(n+1)$ then 
the main case 
(\ref{varch.thm}.1) applies using the  line segment $[(am,0), (i,nj)]$.
The maximum value of $\gamma$ is achieved when $i=0$, giving the condition
$c<a^{-1}+(nd)^{-1}$.

Otherwise $b(i,j)\neq 0$ for some $i+j=md$ and $j< md/(n+1)$
Thus $i<md$. $nj<md$ and  (\ref{varch.thm}.2) applies, giving the condition
$c<d^{-1}$.

Thus either  $(T, \frac{c}{m}(f(x,z^n)=0))$ is klt
as required or we are in the 
 coordinate change case (\ref{varch.thm}.4). 
Thus we consider
$\sum_{ui+vnj=w} b(i,j)x^iz^{nj}$.
Set $U=u/(u,n), V=vn/(u,n), W=w/(u,n)$
and look at the irreducible factorization
$$
\sum_{Ui+Vj=W} b(i,j)x^iy^j=
\prod_{\ell} (\alpha_{\ell} x^V+\beta_{\ell} y^U)^{e_{\ell}}.
$$
Correspondingly, 
$$
\sum_{ui+vnj=w} b(i,j)x^iz^{nj}=
\prod_{\ell} (\alpha_{\ell} x^V+\beta_{\ell} z^{Un})^{e_{\ell}}.
$$
Since (\ref{varch.thm}.4) applies, either $\min\{V,Un\}=1$, that is 
$V=1$, or some $(\alpha_{\ell} x^V+\beta_{\ell} z^{Un})$
factors further. However, in the latter case we have
$(V,Un)>1$ factors of the same multiplicity
$e_{\ell}>w/(v+un)$, which is impossible.
Thus 
$V=1$ and  the coordinate change is
of the form $(x,z)\mapsto (x-\gamma z^{Un}, z)$
which can be realized on $S$ as
$(x,y)\mapsto (x-\gamma y^{U}, y)$.\qed
\medskip

The following consequences show that the
klt condition is satisfied for the cases
leading to $H_2(L,\z)=(\z/m)^2$ for $m\geq 7$
and for $H_2(L,\z)=(\z/3)^6$. This is why I expect 
to have only finitely many families for them.

\begin{cor}\label{klt.cond.cor1} Let $C\subset S$ be a smooth elliptic curve
where $S$ is  $\p^2, \p^1\times \p^1$ or $\f_2$.
Let $D$ be an effective $\q$-divisor such that
 $D\equiv \frac1{n}C$. Then
$(S, (1-\frac1{n})C+D)$ is klt for $n\geq 7$.
\end{cor}

Proof. Start with the $S=\p^2$ case.
If $H\subset \p^2$ is a line then
$C\equiv 3H$ and so $D\equiv \frac3{n}H$.
Thus $(C\cdot D)_p\leq \frac9{n}$ and 
$\mult_pD\leq (H\cdot D)=\frac3{n}$. Since
$$
\tfrac{n}9+\tfrac13>1\qtq{for $n\geq 7$,}
$$
the rest   follows from (\ref{klt.crit}).
The other two cases are similarly easy.
\qed

\begin{cor} Let $C\subset \f_4$ be a smooth curve
of genus 3 such that $C\equiv 2E+8F$.
Let $D$ be an effective $\q$-divisor such that
 $D\equiv -(K+\frac23C)$. Then
$(\f_3, \frac23C+D)$ is klt.
\end{cor}

Proof. Here  $D\equiv \frac23(E+F)$, 
thus 
$$
(C\cdot D)_p\leq \tfrac23\cdot (2E+8F)\cdot (E+F)=\tfrac43,
$$ 
and 
$\mult_pD\leq \frac23$. Since
$\tfrac34+\frac12>1$
the rest   follows from (\ref{klt.crit}).
\qed
\medskip

One can also use (\ref{klt.crit}) to check
 the existence condition \cite{dk} for orbifold
K\"ahler--Einstein metrics.

\begin{cor}\label{ke.conds} Let $S$ be a Del Pezzo surface
with quotent singularities and $C\in |-K_S|$ a smooth elliptic curve.
Then $(S, (1-\frac1{n})C)$ has an orbifold
K\"ahler--Einstein metric whenever $n>\frac23 K_S^2$.
\end{cor}

Proof. Set $d=K_S^2$. Let $D$ be an effective $\q$-divisor such that
 $D\equiv \frac1{n}C$.  We need to check that
$(S, (1-\frac1{n})C+\frac23D)$ is klt.
As in (\ref{klt.cond.cor1}), this holds if
$$
\tfrac23 < \min\left\{\tfrac{n}{d}+\tfrac{n}{nd}, \tfrac{n}{d}\right\}=
\tfrac{n}{d}.\qed
$$

 \begin{ack}  I thank Ch.\ Boyer and K.\ Galicki for many comments
and corrections and the University of Utah where some of this work was done.
Partial financial support  was provided by  the NSF under grant number 
DMS-0500198. 
\end{ack}

\bibliography{refs}

\vskip1cm

\noindent Princeton University, Princeton NJ 08544-1000

\begin{verbatim}kollar@math.princeton.edu\end{verbatim}

\end{document}